\documentclass[12pt,english]{smfart}

\usepackage{smfthm}
\usepackage{amssymb}
\usepackage{mathrsfs}

\newcommand{\ZZ}{\mathbf{Z}}
\newcommand{\FF}{\mathbf{F}}
\newcommand{\RR}{\mathbf{R}}
\newcommand{\QQ}{\mathbf{Q}}

\newcommand{\Qp}{\mathbf{Q}_p}
\newcommand{\Zp}{\mathbf{Z}_p}

\newcommand{\OO}{\mathcal{O}}
\newcommand{\MM}{\mathfrak{m}}

\newcommand{\Qbar}{\overline{\mathbf{Q}}}
\newcommand{\Qpbar}{\overline{\mathbf{Q}}_p}
\newcommand{\Zpbar}{\overline{\mathbf{Z}}_p}
\newcommand{\Qpnr}{\mathbf{Q}_p^{\nr}}
\newcommand{\Qpnrhat}{\widehat{\mathbf{Q}}_p^{\nr}}
\newcommand{\Fbar}{\overline{\mathbf{F}}}
\newcommand{\Fpbar}{\overline{\mathbf{F}}_p}

\newcommand{\eps}{\varepsilon}
\renewcommand{\phi}{\varphi}
\newcommand{\rhobar}{\overline{\rho}}
\renewcommand{\geq}{\geqslant}
\renewcommand{\leq}{\leqslant} 

\newcommand{\Gal}{\operatorname{Gal}}
\newcommand{\galp}{\Gal(\Qpbar/\Qp)}
\newcommand{\galk}{\Gal(\Qpbar/K)}

\newcommand{\frob}{\operatorname{Frob}}
\newcommand{\Mat}{\operatorname{Mat}}
\newcommand{\B}{\operatorname{B}}
\newcommand{\GL}{\operatorname{GL}}

\newcommand{\vp}{\operatorname{val}_p}

\newcommand{\smat}[1]{\left( \begin{smallmatrix} #1 \end{smallmatrix} \right)}
\newcommand{\pmat}[1]{\begin{pmatrix} #1 \end{pmatrix}}
\newcommand{\dpar}[1]{(\!( #1 )\!)}
\newcommand{\dcroc}[1]{[\![ #1 ]\!]}

\newcommand{\bdr}{\mathbf{B}_{\dR}}  
\newcommand{\bst}{\mathbf{B}_{\st}} 
\newcommand{\bcris}{\mathbf{B}_{\cris}} 
\newcommand{\btrig}{\widetilde{\mathbf{B}}^+_{\rig}} 
\newcommand{\btdagrig}{\widetilde{\mathbf{B}}^\dagger_{\rig}} 

\newcommand{\dcris}{\mathrm{D}_{\cris}}
\newcommand{\dst}{\mathrm{D}_{\st}}

\newcommand{\ddr}{\mathrm{D}_{\dR}}

\newcommand{\calE}{\mathcal{E}}
\newcommand{\calR}{\mathcal{R}}
\newcommand{\calS}{\mathcal{S}}
\newcommand{\calC}{\mathcal{C}}
\newcommand{\Gm}{\mathbf{G}_\mathrm{m}}
\newcommand{\bfont}{\mathbf{B}}
\newcommand{\dfont}{\mathrm{D}}

\newcommand{\Hco}{\mathrm{H}}
\newcommand{\Proj}{\mathbf{P}}
\newcommand{\Linv}{\mathcal{L}}
\newcommand{\scrX}{\mathscr{X}}
\newcommand{\scrW}{\mathscr{W}}
\newcommand{\scrU}{\mathscr{U}}

\newcommand{\Fil}{\operatorname{Fil}}
\newcommand{\rk}{\operatorname{rk}}

\newcommand{\End}{\operatorname{End}}
\newcommand{\rig}{\operatorname{rig}}
\newcommand{\irr}{\operatorname{irr}}
\newcommand{\nr}{\operatorname{nr}}
\newcommand{\ab}{\operatorname{ab}}
\newcommand{\cris}{\operatorname{cris}}
\newcommand{\st}{\operatorname{st}}
\newcommand{\et}{\operatorname{\acute{e}t}}

\newcommand{\dR}{\operatorname{dR}}
\newcommand{\cycl}{\operatorname{cycl}}

\newcommand{\ngeo}{\operatorname{ng}}

\newcommand{\Ext}{\operatorname{Ext}}
\newcommand{\Card}{\operatorname{Card}}

\title{Trianguline representations}

\author{Laurent Berger}

\address{UMPA ENS de Lyon \\
UMR 5669 du CNRS \\
Universit\'e de Lyon}
\email{laurent.berger@ens-lyon.fr}

\subjclass{11-02; 11F11; 11F33; 11F80; 11F85; 11S20; 14G22}

\keywords{Galois representations; $p$-adic Hodge theory; $(\phi,\Gamma)$-modules; trianguline representations; $p$-adic local Langlands correspondence; overconvergent modular forms}

\NumberTheoremsIn{subsection}

\begin{document}

\begin{abstract}
Trianguline representations are a certain class of $p$-adic representations of $\galp$, like the crystalline, semistable and de Rham representations of Fontaine. Their definition involves the theory of $(\varphi,\Gamma)$-modules. In this survey, we explain the theory of $(\varphi,\Gamma)$-modules and the definition and properties of trianguline representations. After that, we give some examples of their occurrence in arithmetic geometry.
\end{abstract}

\maketitle

\setcounter{tocdepth}{2}
\tableofcontents

\setlength{\baselineskip}{18pt}

\section{Introduction}

\subsection{Representations of $\Gal(\Qbar/\QQ)$}
\label{rpgqgl}

The starting point for this survey is that one can attach representations of the group $\Gal(\Qbar/\QQ)$ to some objects which occur in arithmetic geometry, for example elliptic curves and modular forms. Suppose for instance that $A$ is an elliptic curve defined over $\QQ$ and choose a prime number $p$. The group $\Gal(\Qbar/\QQ)$ acts on the $p^n$-th torsion points $A[p^n](\Qbar)$ of $A$ and this gives rise to the Tate module of $A$. Tensoring the Tate module with $\Qp$, we get a $2$-dimensional $\Qp$-vector space $V_pA$ which is the $p$-adic representation of $\Gal(\Qbar/\QQ)$ attached to $A$. 

Let $\ell$ be a prime number and choose an embedding $\iota_\ell : \Qbar \to \Qbar_\ell$. This gives rise to a map $\Gal(\Qbar_\ell/\QQ_\ell) \to \Gal(\Qbar/\QQ)$ which is injective and whose image is the decomposition group $D_\ell$ of a place above $\ell$ (a different choice of $\iota_\ell$ gives rise to another subgroup of $\Gal(\Qbar/\QQ)$ which is conjugate to $D_\ell$). The group $D_\ell$ contains the inertia subgroup $I_\ell$ and the quotient $D_\ell / I_\ell$ is isomorphic to $\Gal(\Fbar_\ell/\FF_\ell)$. The group $\Gal(\Fbar_\ell/\FF_\ell)$ is isomorphic to $\widehat{\ZZ}$ and is topologically generated by the Frobenius map $\frob_\ell = [ z \mapsto z^\ell]$. We then have the following theorem which says that the representation $V_pA$ is also ``attached to $A$'' in a deeper way. 

\begin{theo}\label{grellc}
If $\ell \nmid p \cdot \mathrm{Disc}(A)$, then the restriction of $V_p A$ to $I_\ell$ is trivial and $\det(X-\frob_\ell \mid V_pA) = X^2 - a_\ell X + \ell$, where $a_\ell = \ell+1-\Card(A(\FF_\ell))$.
\end{theo}

As $\ell$ runs through a set of primes of density $1$, the groups $D_\ell$ and their conjugates form a dense subset of $\Gal(\Qbar/\QQ)$ by Chebotarev's theorem and therefore theorem \ref{grellc} determines the semisimplification of $V_pA$. If $\ell \neq p$ but $\ell \mid \mathrm{Disc}(A)$, then we also have a description of $V_p A \mid_{D_\ell}$ which now depends on the geometry of $A \bmod{\ell}$. A much deeper problem is the description of the restriction of  $V_p A$ to $D_p$ and this is one of the goals of Fontaine's theory, which we discuss in this survey. 

Let us recall that one can also attach $p$-adic representations of $\Gal(\Qbar/\QQ)$ to modular forms as follows. Let $f=\sum_{n \geq 0} a_n q^n$ be a normalized modular eigenform of weight $k$, level $N$ and character $\eps$, and let $E$ be the field generated over $\Qp$ by the images of the $a_n$ in $\Qpbar$ under some embedding. The field $E$ is a finite extension of $\Qp$ and we have the following theorem, which combines results of Deligne, Eichler-Shimura and Igusa (see theorem 6.1 of Deligne-Serre \cite{DS74}).

\begin{theo}\label{dsfm}
There exists a semisimple $2$-dimensional $E$-linear representation $V_p f$ of $\Gal(\Qbar/\QQ)$ such that for every prime number $\ell \nmid pN$, the restriction of $V_p f$ to $I_\ell$ is trivial and $\det(X-\frob_\ell \mid V_p f) = X^2 - a_\ell X + \eps(\ell)\ell^{k-1}$.
\end{theo}

If $f$ is of weight $2$, then it corresponds to an elliptic curve $A$ and the representation $V_p f$ is the same as the representation $V_p A$ defined above.

\subsection{Trianguline representations and $(\phi,\Gamma)$-modules}
\label{fontint}

Let $E$ be a finite extension of $\Qp$ which is the field of coefficients of the representations we consider. The goal of Fontaine's theory is to study the $E$-linear representations of $\galp$. These may arise as the restriction to $D_p$ of representations of $\Gal(\Qbar/\QQ)$ as above, but they are also interesting considered on their own. A $p$-adic representation of $\galp$ is then a finite dimensional $E$-vector space $V$, along with a continuous $E$-linear action of $\galp$. 

Fontaine's approach has been to construct some ``rings of periods'', for example $\bcris$, $\bst$ and $\bdr$, and to use these rings to define and study crystalline, semistable and de Rham representations (see \S\ref{fonsec} for reminders about this). These constructions allow one to give a complete description of the restriction to $D_p$ of the representations $V_p A$ and $V_p f$ of \S\ref{rpgqgl} (see \S \ref{ptsec}). Another construction of Fontaine's which is crucial in this survey is the theory of $(\phi,\Gamma)$-modules. There are three variants of this theory, and we now describe (and will describe again in more detail in \S\S \ref{pgsec}--\ref{kedsec}) the theory of $(\phi,\Gamma)$-modules over the Robba ring. 

Let $\calR$ be the Robba ring, that is the ring of power series $f(X)= \sum_{n \in \ZZ} a_n X^n$ where $a_n \in E$ and for which there exists $\rho(f)$ such that $f(X)$ converges on the $p$-adic annulus $\rho(f) \leq |X|_p < 1$. This ring is endowed with a Frobenius $\phi$ given by $(\phi f)(X)=f((1+X)^p-1)$ and with an action of $\Zp^\times$ (now called $\Gamma$) given by $([a] f)(X)=f((1+X)^a-1)$ if $a \in \Zp^\times$.

A $(\phi,\Gamma)$-module is a free $\calR$-module of finite rank $d$, endowed with a semilinear Frobenius $\phi$ such that $\Mat(\phi)$ (the matrix of $\phi$ in some basis) belongs to $\GL_d(\calR)$, and with a commuting semilinear continuous action of $\Gamma$. The main result relating $(\phi,\Gamma)$-modules and $p$-adic Galois representations is the following (it combines theorems of Fontaine, Fontaine-Wintenberger, Cherbonnier-Colmez and Kedlaya). Let $\OO_{\calE}^\dag$ be the set of $f(X) \in \calR$ with $|a_n|_p \leq 1$ for all $n \in \ZZ$. We say that a $(\phi,\Gamma)$-module is \'etale if there exists a basis in which $\Mat(\phi) \in \GL_d(\OO_{\calE}^\dag)$. The ring $\btdagrig$ below denotes one of Fontaine's rings of periods.

\begin{theo}\label{pgint}
If $\dfont$ is an \'etale $(\phi,\Gamma)$-module, then $V(\dfont)=(\btdagrig \otimes_{\calR} \dfont)^{\phi=1}$ is a $p$-adic representation of $\galp$, and the resulting functor $\dfont \mapsto V(\dfont)$ gives rise to an equivalence of categories: \{\'etale $(\phi,\Gamma)$-modules\} $\to$ \{$p$-adic representations\}.
\end{theo}

We denote by $V \mapsto \dfont(V)$ the inverse functor. The category of \'etale $(\phi,\Gamma)$-modules is a full subcategory of the larger category of all $(\phi,\Gamma)$-modules. In particular, if $V$ is an irreducible $p$-adic representation, then $\dfont(V)$ is irreducible in the category of \'etale $(\phi,\Gamma)$-modules but it can be reducible in the larger category of all $(\phi,\Gamma)$-modules.

\begin{defi}\label{dti}
If $V$ is a $p$-adic representation of $\galp$, then we say that $V$ is trianguline if $\dfont(V)$ is a successive extension of $(\phi,\Gamma)$-modules of rank $1$ (after possibly enlarging $E$).
\end{defi}

This definition was first given by Colmez in his construction of the ``unitary principal series of $\GL_2(\Qp)$'', which is an important building block of the $p$-adic local Langlands correspondence for $\GL_2(\Qp)$ (see \S \ref{llsec}). Some important examples of trianguline representations are (1) the semistable representations of $\galp$ and (2) the restriction to $\galp$ of the representations of $\Gal(\Qbar/\QQ)$ attached to finite slope overconvergent modular forms.

This survey has three chapters. In the first one, we give a more detailed description of the definition and properties of $(\phi,\Gamma)$-modules, including Kedlaya's theory of Frobenius slopes. In the second one, we give some examples of trianguline representations by relating the theory of $(\phi,\Gamma)$-modules to $p$-adic Hodge theory, and then we give Colmez' construction of a parameter space for all $2$-dimensional trianguline representations. In the last chapter, we explain how trianguline representations occur in the $p$-adic local Langlands correspondence, in the theory of overconvergent modular forms and in the study of Selmer groups.

\subsection{Notations and conventions}
\label{nc}

The field $E$ is a finite extension of $\Qp$ with ring of integers $\OO_E$ whose maximal ideal is $\MM_E$ and whose residue field is $k_E$. All the characters, representations and group actions in this survey are assumed to be continuous (note that a character $\delta : \Qp^\times \to E^\times$ is necessarily continuous by exercise 6 of \S 4.2 of \cite{S94}). When we say that an $E$-linear object is irreducible, we mean that it is absolutely irreducible, meaning that it remains irreducible when we extend scalars from $E$ to a finite extension.

The cyclotomic character $\chi_{\cycl}$ gives an isomorphism $\chi_{\cycl} : \Gal(\Qp(\mu_{p^\infty}) / \Qp) \to \Zp^\times$. The maximal abelian extension of $\Qp$ is $\Qp^{\ab} = \Qp^{\nr} \cdot \Qp(\mu_{p^\infty})$, and every element of $\Gal(\Qp^{\ab}/\Qp)$ can be written as $\frob_p^n \cdot g$ where $\frob_p$ is the lift of $[z \mapsto z^p]$ and $n \in \widehat{\ZZ}$ and $g \in \Gal(\Qp^{\ab} / \Qp^{\nr})$. If $\delta : \Qp^\times \to \OO_E^\times$ is a unitary character, then by local class field theory $\delta$ gives rise to a character (still denoted by $\delta$) of $\galp$ which is determined by the formula $\delta(\frob_p^n \cdot g) = \delta(p)^{-n} \cdot\delta(\chi(g))$ if $n \in \ZZ$. In other words, we normalize class field theory so that $p$ corresponds to the geometric Frobenius $\frob_p^{-1}$.

\section{Galois representations and $(\phi,\Gamma)$-modules}
\label{pgch}

In this chapter, we explain the theory of $(\phi,\Gamma)$-modules and its relation to $p$-adic representations. This allows us to define trianguline representations.

\subsection{The Robba ring and $(\phi,\Gamma)$-modules}\label{pgsec}

The \emph{Robba ring} $\calR$ is the ring of power series $f(X) = \sum_{n \in \ZZ} a_n X^n$ where $a_n \in E$ such that $f(X)$ converges on an annulus of the form $\rho(f) \leq |X|_p < 1$. For example, the power series $t = \log(1+X)$ belongs the Robba ring (and here one can take $\rho(t)=0$). 

The Robba ring is endowed with a Frobenius map $\phi$ given by $(\phi f)(X) = f((1+X)^p-1)$. Let $\Gamma$ be another notation for $\Zp^\times$ with the isomorphism $\Zp^\times \to \Gamma$ denoted by $a \mapsto [a]$. The Robba ring is endowed with an action of $\Gamma$ given by $([a]f)(X) = f((1+X)^a-1)$ and this action commutes with $\phi$. For example, we have $\phi(t)=pt$ and $[a](t) = at$.

\begin{defi}\label{dfpgr}
A \emph{$(\phi,\Gamma)$-module over $\calR$} is a free $\calR$-module of finite rank $d$, endowed with a semilinear Frobenius $\phi$ such that $\Mat(\phi)$ (the matrix of $\phi$ in some basis) belongs to $\GL_d(\calR)$ and a semilinear action of $\Gamma$ which commutes with $\phi$. 
\end{defi}

There is then an obvious notion of morphism of $(\phi,\Gamma)$-modules, and this gives rise to the category of $(\phi,\Gamma)$-modules. This category is not abelian, since the quotient of two $(\phi,\Gamma)$-modules is not necessarily free (consider for instance the map: $t \cdot \calR \to \calR$).

If $\delta : \Qp^\times \to E^\times$ is a character, then we define $\calR(\delta)$ as the $(\phi,\Gamma)$-module of rank $1$ having $e_\delta$ as a basis where $\phi(e_\delta)=\delta(p)e_\delta$ and $[a](e_\delta) = \delta(a)e_\delta$. The following theorem is proposition 3.1 of \cite{PC08}.

\begin{theo}\label{pgr1}
Every $(\phi,\Gamma)$-module of rank $1$ over $\calR$ is isomorphic to $\calR(\delta)$ for a well-defined character $\delta : \Qp^\times \to E^\times$.
\end{theo}

\subsection{Etale $(\phi,\Gamma)$-modules and Galois representations}\label{rosec}

The ring $\calE^{\dagger}$ is the subring of $\calR$ consisting of those power series $f(X) = \sum_{n \in \ZZ} a_n X^n$ for which the sequence $\{a_n\}_{n \in \ZZ}$ is bounded. The subring of $\calE^\dag$ consisting of those $f(X) = \sum_{n \in \ZZ} a_n X^n$ for which $|a_n|_p \leq 1$ is denoted by $\OO_{\calE}^\dagger$. This is a henselian local ring with residue field $k_E \dpar{X}$.

\begin{defi}\label{dfet}
We say that a $(\phi,\Gamma)$-module over $\calR$ is \emph{\'etale} if it has a basis in which $\Mat(\phi) \in \GL_d(\OO_{\calE}^\dag)$. 
\end{defi}

In \S 2.3 of \cite{LB02}, a ring $\btdagrig$ is constructed which has the following properties: it is endowed with a Frobenius $\phi$ and a commuting action of $\galp$, and it contains the Robba ring $\calR$. This inclusion is compatible with $\phi$ and with the action of $\Gamma$ on $\calR$, in the sense that if $y \in \calR$ and $g \in \galp$, then $g(y) = [\chi_{\cycl}(g)](y)$. One can think of $\btdagrig$ as some sort of ``algebraic closure'' of $\calR$.

If $\dfont$ is a $(\phi,\Gamma)$-module over $\calR$, then $V(\dfont) = (\btdagrig \otimes_{\calR} \dfont)^{\phi=1}$ is an $E$-vector space, endowed with the action of $\Gal(\Qpbar/\Qp)$ given by $g(x \otimes e) = g(x) \otimes [\chi_{\cycl}(g)](e)$. This $E$-vector space can be finite or infinite-dimensional in general, but we have the following theorem which combines results of Fontaine (theorem 3.4.3 of \cite{F90}), Cherbonnier-Colmez (corollary III.5.2 of \cite{CC98}) and Kedlaya (theorem 6.3.3 of \cite{KK05}).

\begin{theo}\label{pgrep}
If $\dfont$ is an \'etale $(\phi,\Gamma)$-module of rank $d$ over $\calR$, then $V(\dfont)$ is an $E$-linear representation of dimension $d$ of $\Gal(\Qpbar/\Qp)$,  and the resulting functor, from the category of \'etale $(\phi,\Gamma)$-modules over $\calR$ to the category of $E$-linear representations of $\Gal(\Qpbar/\Qp)$, is an equivalence of categories.
\end{theo}

We denote by $V \mapsto \dfont(V)$ the inverse functor, which to a $p$-adic representation attaches the corresponding \'etale $(\phi,\Gamma)$-module over $\calR$.

For example, the $(\phi,\Gamma)$-module $\calR(\delta)$ is \'etale if and only if $\vp(\delta(p))=0$. In this case, the representation $V(\calR(\delta))$ is the character of $\galp$ corresponding to $\delta$ by local class field theory, as recalled in \S\ref{nc}.

\subsection{Trianguline representations}\label{trsec} 

We can now give the definition of trianguline representations (see \S 0.4 of \cite{PC08}).

\begin{defi}\label{dtr}
If $V$ is a $p$-adic representation of $\galp$, then 
\begin{enumerate}
\item we say that $V$ is \emph{split trianguline} if the $(\phi,\Gamma)$-module $\dfont(V)$ is a successive extension of $(\phi,\Gamma)$-modules of rank $1$;
\item we say that $V$ is \emph{trianguline} if there exists a finite extension $F$ of $E$ such that $F \otimes_E V$ is split trianguline.
\end{enumerate}
\end{defi}

In other words, a $p$-adic representation $V$ is split trianguline if and only if $\dfont(V)$ has a basis in which the matrices of $\phi$ and of the elements of $\Gamma$ are all upper-triangular.

On the level of $(\phi,\Gamma)$-modules, the possible extension of scalars from $E$ to $F$ consists in extending scalars from the Robba ring with coefficients in $E$ to the Robba ring with coefficients in $F$. For example, we'll see later on that semistable representations are always trianguline, and they are split trianguline if and only if $E$ contains the eigenvalues of $\phi$ on $\dst(V)$.

It is important to understand that a representation $V$ may well be trianguline without $V$ itself being an extension of representations of dimension $1$. Indeed, the definition is that $\dfont(V)$ is a successive extension of $(\phi,\Gamma)$-modules of rank $1$, but these $(\phi,\Gamma)$-modules are generally not \'etale and therefore do not correspond to subquotients of $V$.

Note also that a $(\phi,\Gamma)$-module may be written as a successive extension of $(\phi,\Gamma)$-modules of rank $1$ in several different ways.

In the rest of this survey, we'll see several examples of trianguline representations, but we give here the two main classes:

\begin{enumerate}
\item the representations of $\galp$ which become semistable when restricted to $\Gal(\Qpbar/\Qp(\zeta_{p^n}))$ for some $n \geq 0$;
\item the representations of $\galp$ which arise from overconvergent modular eigenforms of finite slope.
\end{enumerate}

In \cite{LBDP}, some explicit families of $2$-dimensional representations are constructed and the trianguline ones are determined.

\subsection{Slopes of $(\phi,\Gamma)$-modules}\label{kedsec}

We now recall Kedlaya's theory of slopes for $\phi$-modules over the ring $\calR$ (i.e.\ free $\calR$-modules of finite rank $d$ with a semilinear $\phi$ such that $\Mat(\phi) \in \GL_d(\calR)$). If $s = a/h \in \QQ$ is written in lowest terms, then we say that a $\phi$-module over $\calR$ is \emph{pure of slope $s$} if it is of rank $\geq 1$ and has a basis in which $\Mat(p^{-a} \phi^h) \in \GL_d(\OO_{\calE}^\dag)$ (being \'etale is therefore equivalent to being pure of slope zero). A $\phi$-module over $\calR$ which is pure of a certain slope is said to be \emph{isoclinic}. For example, the $(\phi,\Gamma)$-module $\calR(\delta)$ is pure of slope $\vp(\delta(p))$. The main result of the theory of slopes is theorem 6.10 of \cite{KK04}.

\begin{theo}\label{kedsl}
If $\dfont$ is a $\phi$-module over $\calR$, then there exists a unique filtration $\{0\} = \dfont_0 \subset \dfont_1 \subset \cdots \subset \dfont_\ell = \dfont$ of $\dfont$ by sub-$\phi$-modules such that:
\begin{enumerate}
\item for all $i \geq 1$, $\dfont_i / \dfont_{i-1}$ is an isoclinic $\phi$-module;
\item if $s_i$ is the slope of $\dfont_i / \dfont_{i-1}$, then $s_1 < s_2 < \cdots < s_\ell$.
\end{enumerate}
\end{theo}

If $\dfont$ is a $(\phi,\Gamma)$-module, then each of the $\dfont_i$ is stable under the action of $\Gamma$ since the filtration is unique, and hence each $\dfont_i$ is itself a $(\phi,\Gamma)$-module. 

A delicate but crucial point of the theory of slopes is that a $\phi$-module over $\calR$ which is pure of slope $s$ has no subobject of slope $<s$ by theorem \ref{kedsl}, but it may well have  subobjects of slope $> s$. This helps to explain the definition of trianguline representations: an \'etale $(\phi,\Gamma)$-module over $\calR$ may be irreducible in the category of \'etale $(\phi,\Gamma)$-modules but it can still admit some nontrivial subobjects in the larger category of all $(\phi,\Gamma)$-modules.

Theorem \ref{kedsl} also helps to understand theorem \ref{pgrep}. If $\dfont$ is a $(\phi,\Gamma)$-module, then $V(\dfont) = (\btdagrig \otimes_{\calR} \dfont)^{\phi=1}$ is constructed by solving $\phi$-equations determined by the matrix of $\phi$ on $\dfont$. If the slopes of $\dfont$ are $>0$, then these equations have no nonzero solutions, while if the slopes of $\dfont$ are $<0$, then the space of solutions if infinite dimensional (see theorem A of \cite{LB09} for more precise results). The condition that $\dfont$ is \'etale is exactly the right one for $V(\dfont)$ to be a finite dimensional $E$-vector space of the correct dimension.

\section{Examples of trianguline representations}
\label{exch}

In this chapter, we explain how to relate $(\phi,\Gamma)$-modules and $p$-adic Hodge theory, which allows us to give important examples of trianguline representations. After that, we explain how to compute extensions of $(\phi,\Gamma)$-modules and Colmez' resulting construction of all $2$-dimensional trianguline representations.

\subsection{Fontaine's rings of periods}\label{fonsec}

The purpose of Fontaine's theory is to sort through $p$-adic representations, and to classify the interesting ones by using objects from linear algebra. Recall that Fontaine has constructed in \cite{FP} a number of rings, for example $\bcris$, $\bst$ and $\bdr$. The construction of these rings is quite complicated but they have a number of properties, some of which we now recall and which suffice for this survey. All of them are $\Qp$-algebras endowed with an action of $\galp$ and some extra structures, which are all compatible with the action of $\galp$. The ring $\bst$ has a Frobenius $\phi$ and a monodromy operator $N$ which satisfy the relation $N \phi = p \phi N$, and the ring $\bcris$ is then $\bst^{N=0}$. The ring $\bdr$ is actually a field, and is endowed with a filtration. The ring $\bcris$ contains $\Qpnrhat$ and the choice of $\log_p(p)$ gives rise to an injective map $\Qpbar \otimes_{\Qpnr} \bst \to \bdr$.

If $V$ is a $p$-adic representation of $\galp$ and $\ast \in \{$cris, st, dR$\}$, then we set $\dfont_\ast(V) = (\bfont_\ast \otimes_{\Qp} V)^{\galp}$. The space $\dfont_\ast(V)$ is then an $E$-vector space of dimension $\leq \dim_E(V)$ and we say that $V$ is \emph{crystalline} or \emph{semistable} or \emph{de Rham} if we have equality of dimensions with $\ast$ being cris, st or dR.

The $E$-vector space $\ddr(V)$ is then endowed with an $E$-linear filtration, the space $\dst(V) \subset \ddr(V)$ is a filtered $(\phi,N)$-module (that is, a finite dimensional $E$-vector space with an $E$-linear map $\phi$, an $E$-linear map $N$ satisfying the relation $N \phi = p \phi N$, and a filtration by $E$-vector subspaces, which is not assumed to be stable under either $\phi$ or $N$) and $\dcris(V)=\dst(V)^{N=0}$ is a filtered $\phi$-module. 

If $D$ is a filtered $(\phi,N)$-module, then we define the Newton number $t_N(D)$ as the $p$-adic valuation of $\phi$ on $\det(D)$ and the Hodge number $t_H(D)$ as the unique integer $h$ such that $\Fil^h(\det(D))=\det(D)$ and $\Fil^{h+1}(\det(D))=\{0\}$. We say that $D$ is \emph{admissible} if $t_H(D)=t_N(D)$ and if $t_H(D') \leq t_N(D')$ for every $(\phi,N)$-stable subspace $D'$ of $D$. The following theorem combines results of Fontaine (\S 5.4 of \cite{FST}) and the Colmez-Fontaine theorem (theorem A of \cite{CF00}).

\begin{theo}\label{cf}
If $V$ is a semistable representation of $\galp$, then $\dst(V)$ is an admissible filtered $(\phi,N)$-module, and the functor $V \mapsto \dst(V)$ gives an equivalence of categories: \{semistable representations\} $\to$ \{admissible filtered $(\phi,N)$-modules\}.
\end{theo}

All of these constructions also work for representations of $\galk$ if $K$ is a finite extension of $\Qp$. In particular, we say that a $p$-adic representation of $\galp$ is \emph{potentially semistable} if its restriction to $\galk$ is semistable for some finite extension $K$ of $\Qp$. Potentially semistable representations of $\galp$ are always de Rham.

\subsection{$p$-adic Hodge theory}\label{ptsec}

If $X$ is a proper and smooth scheme over $\Qp$, then the \'etale cohomology groups $\Hco^i_{\et}(X_{\Qpbar},\Qp)$ are $p$-adic representations of $\galp$. If $X$ has good reduction at $p$, then we can consider its crystalline cohomology groups which have the structure of filtered $\phi$-modules, and if $X$ has bad semistable reduction at $p$, then one can replace the crystalline cohomology groups with a generalization: the log-crystalline cohomology groups, which have the structure of filtered $(\phi,N)$-modules. We then have the following theorem of Tsuji (theorem 0.2 of \cite{TT99}), which is the former conjecture $C_{\st}$ of Fontaine-Jannsen (see \S 6.2 of \cite{FST}).

\begin{theo}\label{tsu}
If $X$ is a proper scheme over $\Zp$ with semistable reduction, then $\Hco^i_{\mathrm{\et}}(X_{\Qpbar},\Qp)$ is a semistable representation of $\galp$, and there is a natural isomorphism of filtered $(\phi,N)$-modules: $\dst(\Hco^i_{\mathrm{\et}}(X_{\Qpbar},\Qp)) = \Hco^i_{\mathrm{log}\text{-}\cris}(X)$.
\end{theo}

If $f$ is a modular eigenform, then one can attach to it a $p$-adic representation $V_p f$ as recalled in theorem \ref{dsfm}. The representation $V_p f$ is always potentially semistable and a result of Saito (the main theorem of \cite{TS97}) completely describes the restriction of $V_p f$ to $D_p$. If $p \nmid N$, then Saito's theorem was previously proved by Scholl (see theorem 1.2.4 of \cite{AS90}). In this case, $V_p f$ is crystalline and $\dcris((V_p f)^*) = D_{k,a_p}$ where $k$ is the weight of $f$, $a_p$ is the eigenvalue of the Hecke operator $T_p$, and $D_{k,a_p} = E e_1 \oplus E e_2$ with
\[ \Mat(\phi) = \pmat{0 & -1 \\ \eps(p)p^{k-1} & a_p} 
\text{ and } \Fil^i D_{k,a_p} = 
\begin{cases} 
D_{k,a_p} & \text{ if $i \leq 0$,} \\
E e_1 & \text{ if $1 \leq i \leq k-1$,} \\
\{0\} & \text{ if $i \geq k$.}
\end{cases} \]

The following is known as the Fontaine-Mazur conjecture (conjecture 1 of \cite{FM95}). 

\begin{conj}\label{fontmaz}
If $V$ is an irreducible $p$-adic representation of $\Gal(\Qbar/\QQ)$, whose restriction to $I_\ell$ is trivial for all $\ell$ except a finite number, and whose restriction to $D_p$ is potentially semistable, then $V$ is a subquotient of an \'etale cohomology group of some algebraic variety over $\QQ$.
\end{conj}

If in addition $\dim(V)=2$ and $V$ is odd, then we actually expect $V$ to be the representation attached to a modular eigenform, and we have the following precise conjecture (conjecture 3c of \cite{FM95}). The \emph{Hodge-Tate weights} of a de Rham representation $V$ are the opposites of the jumps of the filtration on $\ddr(V)$.

\begin{conj}\label{fm}
If $V$ is an irreducible $2$-dimensional $p$-adic representation of $\Gal(\Qbar/\QQ)$, whose restriction to $I_\ell$ is trivial for all $\ell$ except a finite number, and whose restriction to $D_p$ is potentially semistable with distinct Hodge-Tate weights, then $V$ is a twist of the Galois representation attached to a cuspidal eigenform with weight $k \geq 2$.
\end{conj}

Let us write $\overline{V}$ for the reduction modulo $\MM_E$ of $V$.

\begin{theo}\label{emkis}
Conjecture \ref{fm} is true, if we suppose that $\overline{V}$ satisfies some technical hypotheses.
\end{theo}

This theorem has been proved independently by Kisin (it is the main theorem of \cite{MK09}) and by Emerton (theorem 1.2.4 of \cite{ME10}). The ``technical hypotheses'' of Kisin are the following ($\chi_{\cycl}$ is now the reduction mod $p$ of the cyclotomic character, and $*$ denotes a cocycle which may be equal to $0$). 
\begin{enumerate}
\item $p \neq 2$ and $\overline{V}$ is odd,
\item $\overline{V} \mid_{\Gal(\Qbar/\QQ(\zeta_p))}$ is irreducible,
\item $\overline{V} \mid_{\Gal(\Qpbar/\Qp)}$ is not of the form $\smat{\eta \chi_{\cycl} & * \\ 0 & \eta}$ for any character $\eta$.
\end{enumerate}

The ``technical hypotheses'' of Emerton are (1) and (2) and
\begin{itemize}
\item[3'.] $\overline{V} \mid_{\Gal(\Qpbar/\Qp)}$ is not of the form $\smat{\eta & * \\ 0 & \eta \chi_{\cycl}}$ nor of the form $\smat{\eta & * \\ 0 & \eta}$ for any character $\eta$.
\end{itemize}

\subsection{Crystalline and semistable $(\phi,\Gamma)$-modules}\label{bersec}

In \S \ref{fonsec}, we recalled the definition of $\dcris(V)$ and $\dst(V)$ for a $p$-adic representation $V$. We now explain how to extend this definition to $(\phi,\Gamma)$-modules. Recall that we denote by $t$ the element $\log(1+X) \in \calR$. 

\begin{defi}\label{dfdc}
If $\dfont$ is a $(\phi,\Gamma)$-module, let $\dcris(\dfont) = (\calR[1/t] \otimes_{\calR} \dfont)^{\Gamma}$.
\end{defi}

In order to define $\dst(\dfont)$, we add a variable to $\calR$ as follows. The power series $\log(\phi(X)/X^p)$ and $\log(\gamma(X)/X)$ (for $\gamma \in \Gamma$) both converge in $\calR$. Let $\log(X)$ be a variable which we adjoin to $\calR$, with the Frobenius and the action of $\Gamma$ extending to $\calR[\log(X)]$ by $\phi(\log(X)) = p \log(X) + \log(\phi(X)/X^p)$ and $\gamma(\log(X)) = \log(X) + \log(\gamma(X)/X)$. We also define a monodromy map $N$ on $\calR[\log(X)]$ by $N=-p/(p-1) \cdot d/d \log(X)$. 

\begin{defi}\label{dfds}
If $\dfont$ is a $(\phi,\Gamma)$-module, let $\dst(\dfont) = (\calR[\log(X),1/t] \otimes_{\calR} \dfont)^{\Gamma}$.
\end{defi}

Definitions \ref{dfdc} and \ref{dfds} make sense for any $(\phi,\Gamma)$-module. We say that $\dfont$ is \emph{crystalline} or \emph{semistable} if $\dcris(\dfont)$ or $\dst(\dfont)$ is an $E$-vector space of dimension $\rk(\dfont)$. The space $\dst(\dfont)$ is then a $(\phi,N)$-module and $\dcris(\dfont) = \dst(\dfont)^{N=0}$. One can also define a filtration on these two spaces by using the filtration of $\calR$ given by ``the order of vanishing at $\zeta_{p^n}-1$ for $n \gg 0$'' so that $\dst(\dfont)$ becomes a filtered $(\phi,N)$-module (which in general is not admissible). The following result is theorem 0.2 of \cite{LB02}.

\begin{theo}\label{recip}
If $V$ is a $p$-adic representation of $\galp$, and if $\dfont(V)$ is the attached $(\phi,\Gamma)$-module, then $\dcris(V)=\dcris(\dfont(V))$ and $\dst(V)=\dst(\dfont(V))$.
\end{theo}

The proof of this requires a number of delicate computations in several of Fontaine's rings of periods. Recall that $\btdagrig$ is the ring used in \S \ref{rosec} in order to attach $p$-adic representations to $(\phi,\Gamma)$-modules. One can show that the ring $\bcris$ of Fontaine admits a subring $\btrig$ such that 
\begin{enumerate}
\item for any $p$-adic representation $V$, the inclusion $(\btrig[1/t] \otimes_{\Qp} V)^{\galp} \subset \dcris(V)$ is an isomorphism;
\item there is a natural inclusion $\btrig \subset \btdagrig$.
\end{enumerate}
These facts allow one to go from the usual $p$-adic periods to the theory of $(\phi,\Gamma)$-modules, and then to prove theorem \ref{recip}. The spaces $\dst(V)$ and $\dst(\dfont(V))$ are then equal as subspaces of $\btdagrig[1/t] \otimes_{\Qp} V$. It is also possible to define $\ddr(\dfont)$ as well as de Rham $(\phi,\Gamma)$-modules in the same way, and to prove an analogue of theorem \ref{recip}, but this is slightly more complicated and we do not give the recipe here.

If $V$ is a semistable representation, and if $M$ is a $(\phi,N)$-stable subspace of $\dst(V)$, then it is easy to see that $( \calR[\log(X),1/t] \otimes_E M )^{N=0} \cap \dfont(V)$ is a sub $(\phi,\Gamma)$-module of $\dfont(V)$ of rank $\dim(M)$. Using this observation and theorem \ref{recip}, we get the following result.

\begin{theo}\label{sstrig}
Semistable representations of $\galp$ are trianguline.
\end{theo}

We see that the $(\phi,\Gamma)$-module of a semistable representation may then admit several different triangulations, corresponding to flags of $\dst(V)$ stable under $\phi$ and $N$. Another consequence of theorem \ref{recip}, which is proved in the same way, is the following useful result (proposition 4.3 of \cite{PC08}).

\begin{theo}\label{dcnz}
If $V$ is a $p$-adic representation of dimension $2$, then $V$ is trianguline if and only if there exists a character $\eta$ of $\galp$ such that $\dcris(V(\eta)) \neq 0$.
\end{theo}

\subsection{Weights of trianguline representations}\label{wtsec}

Recall that $p$-adic representations of $\galp$ have weights: Sen's theory (\S 2.2 of \cite{SS80}) allows us to attach to $V$ a polynomial $P(X) \in E[X]$ of degree $\dim(V)$, whose roots are the \emph{generalized Hodge-Tate weights of $V$} (warning: in \cite{BC09} as in other places, the opposite sign is chosen for the weights. For us, the cyclotomic character has weight $+1$). For example if $V$ is de Rham, then these weights are the opposites of the jumps of the filtration on $\ddr(V)$, and are then the classical Hodge-Tate weights of $V$. 

If $V$ is a trianguline representation and if $\{0\} = \dfont_0 \subset \dfont_1 \subset \cdots \subset \dfont_d = \dfont(V)$ is a triangulation of $V$, then each $\dfont_i / \dfont_{i-1}$ is of rank $1$ and hence of the form $\calR(\delta_i)$ by theorem \ref{pgr1}. If $\delta : \Qp^\times \to E^\times$ is a character, then $w(\delta)=\log_p \delta(u)/\log_p u$ does not depend on the choice of $u \in 1+p\Zp \setminus \{1\}$ and is called the \emph{weight} of $\delta$. 

\begin{theo}\label{wtri}
If $V$ is a trianguline representation and $\delta_1$, \dots, $\delta_d$ are as above, then $w(\delta_1)$, \dots, $w(\delta_d)$ are the generalized Hodge-Tate weights of $V$.
\end{theo}

The following theorem (proposition 2.3.4 of \cite{BC09}) can be seen as a generalization of Perrin-Riou's theorem 1.5 of \cite{PR94} that ``ordinary representations are semistable''.

\begin{theo}\label{trdr}
Let $V$ be a trianguline representation. If $V$ admits a triangulation with characters $\delta_1$, \dots, $\delta_d$ such that $w(\delta_1)$, \dots, $w(\delta_d)$ are integers and $w(\delta_1) > \cdots > w(\delta_d)$, then $V$ is de Rham.
\end{theo}

\subsection{Cohomology of $(\phi,\Gamma)$-modules}\label{liusec}

Since trianguline representations are successive extensions of $(\phi,\Gamma)$-modules of rank $1$, an important part of the study of these representations is the determination of the extension groups of $(\phi,\Gamma)$-modules.

Let $\dfont$ be a $(\phi,\Gamma)$-module and let $\gamma$ be a topological generator of $\Gamma$ (the group $\Zp^\times$ is topologically cyclic if $p\neq 2$; if $p=2$, then the definitions have to be slightly modified). Let $C(\phi,\gamma)$ be the complex (first considered by Herr in \cite{LH98})
\[ 0 \to \dfont \xrightarrow{z \mapsto ((\gamma-1)z,(\phi-1)z)} \dfont \oplus \dfont \xrightarrow{(x,y) \mapsto (\phi-1)x - (\gamma-1)y} \dfont \to 0. \]
The $E$-vector spaces $\Hco^i(C(\phi,\gamma))$ do not depend on the choice of $\gamma$ and we define the cohomology groups of $\dfont$ to be $\Hco^i(\dfont) = \Hco^i(C(\phi,\gamma))$. Note that by construction $\Hco^i(\dfont)=0$ if $i \geq 3$. 

The following theorem (theorems 1.1 and 1.2 and \S 3.1  of \cite{RL08}) summarizes several properties of the groups $\Hco^i(\dfont)$; we write $h^i(\dfont)$ for $\dim_E \Hco^i(\dfont)$.

\begin{theo}\label{pgcoh}
If $\dfont$ is a $(\phi,\Gamma)$-module, then:
\begin{enumerate}
\item the $\Hco^i(\dfont)$ are finite dimensional $E$-vector spaces and $h^0(\dfont)-h^1(\dfont)+h^2(\dfont)=-\rk(\dfont)$;
\item $\Hco^0(\dfont)=\dfont^{\Gamma=1,\phi=1}$ and $\Hco^1(\dfont)=\Ext^1(\calR,\dfont)$;
\item if $V$ is a $p$-adic representation, then $\Hco^i(\dfont(V)) \simeq \Hco^i(\galp,V)$;
\end{enumerate}
\end{theo}

Combining (3) and (1), we recover Tate's Euler characteristic formula.

In the special case when $\dfont$ is of rank $1$, Colmez has computed explicitly $\Hco^1(\dfont)$. We have the following result (theorem 0.2 of \cite{PC08}) which we use in \S\ref{exsec}. Let $x : \Qp^\times \to E^\times$ be the map $z \mapsto z$ and let $|\cdot|_p : \Qp^\times \to E^\times$ be the map $z \mapsto p^{-\vp(z)}$. 

\begin{theo}\label{extpgr}
If $\delta_1$ and $\delta_2 : \Qp^\times \to E^\times$ are two characters, then $\Ext^1(\calR(\delta_2), \calR(\delta_1))$ is a $1$-dimensional $E$-vector space, unless $\delta_1 \delta_2^{-1}$ is either of the form $x^{-i}$ with $i \geq 0$ or $|x|_px^i$ with $i \geq 1$, in which case $\Ext^1(\calR(\delta_2), \calR(\delta_1))$ is of dimension $2$.
\end{theo}

In the first case, there is therefore one nonsplit extension $0 \to \calR(\delta_1) \to \dfont \to \calR(\delta_2) \to 0$ while in the second case, the set of such extensions is parameterized by $\Proj^1(E)$. The parameter for such an extension is called the \emph{$\Linv$-invariant} and turns out to be a generalization of the usual $\Linv$-invariants (see \cite{PCLI} for a discussion of $\Linv$-invariants of modular forms).

\subsection{Trianguline representations of dimension $2$}\label{exsec}

We now explain how we can use the results of the preceding paragraph in order to construct a parameter space for all irreducible trianguline representations of dimension $2$.

If $\delta : \Qp^\times \to E^\times$ is a character, then we set $u(\delta)=\vp(\delta(p))$, so that $u(\delta)$ is the slope of $\calR(\delta)$ as in \S\ref{kedsec}. Recall that $w(\delta)$ is the weight of $\delta$, defined in \S \ref{wtsec}.

If $V$ is a trianguline representation of dimension $2$, then $\dfont(V)$ is an extension of two $(\phi,\Gamma)$-modules of rank $1$, so that we have an exact sequence $0 \to \calR(\delta_1) \to \dfont(V) \to \calR(\delta_2) \to 0$. The fact that $\dfont(V)$ is \'etale implies that $u(\delta_1) + u(\delta_2)=0$, and (because of theorem \ref{kedsl}) $u(\delta_1) \geq 0$. If $u(\delta_1)=u(\delta_2)=0$, then $\calR(\delta_1)$ and $\calR(\delta_2)$ are \'etale, and $V$ itself is an extension of two representations. 

Let $\calS$ be the set $\{ (\delta_1,\delta_2,\Linv) \}$ where $\delta_1$ and $\delta_2$ are characters $\Qp^\times \to E^\times$, and $\Linv = \infty$ if $\delta_1 \delta_2^{-1}$ is neither of the form $x^{-i}$ with $i \geq 0$ nor of the form $|x|_p x^i$ with $i \geq 1$, and $\Linv \in \Proj^1(E)$ otherwise. Theorem \ref{extpgr} above allows us to construct for every $s \in \calS$ a nontrivial extension $\dfont(s)$ of $\calR(\delta_2)$ by $\calR(\delta_1)$, and vice versa.

If $s \in \calS$, then we set $w(s)=w(\delta_1)-w(\delta_2)$. We define $\calS_*$ as the set of $s \in \calS$ such that $u(\delta_1)+u(\delta_2)=0$ and $u(\delta_1)>0$, and we then set $u(s)=u(\delta_1)$ if $s \in \calS_*$. We define the ``crystalline'', ``semistable'' and ``nongeometric'' parameter spaces as follows.
\begin{enumerate}
\item $\calS^{\cris}_* = \{ s \in \calS_*$ such that $w(s) \in \ZZ_{\geq 1}$ and $u(s)<w(s)$ and $\Linv=\infty \}$;
\item $\calS^{\st}_* = \{ s \in \calS_*$ such that $w(s) \in \ZZ_{\geq 1}$ and $u(s)<w(s)$ and $\Linv \neq \infty \}$;
\item $\calS^{\ngeo}_* = \{ s \in \calS_*$ such that $w(s) \notin \ZZ_{\geq 1} \}$.
\end{enumerate}
Let $\calS_{\irr} = \calS^{\cris}_* \sqcup \calS^{\st}_* \sqcup \calS^{\ngeo}_*$.

\begin{theo}\label{dsirr}
If $s \in \calS_{\irr}$, then $\dfont(s)$ is \'etale, and the attached representation $V(s)$ is trianguline and irreducible. Every $2$-dimensional irreducible trianguline representation is of the form $V(s)$ for some $s \in \calS_{\irr}$ (after possibly extending scalars), and we have $V(s)=V(s')$ if and only if $s \in \calS^{\cris}_*$ and $s'=(x^{w(s)} \delta_2, x^{-w(s)} \delta_1,\infty)$. 
\end{theo}

In particular, if $s \in \calS \setminus \calS_{\irr}$, then either $\dfont(s)$ is \'etale but $V(s)$ is reducible, or $\dfont(s)$ is not even \'etale (this happens for example if $w(s) \in \ZZ_{\geq 1}$ and $u(s)>w(s)$). These cases are examined in \S 3 of \cite{PC08}.

If $s \in \calS^{\cris}_*$, then the representation $V(s)$ becomes crystalline over an abelian extension of $\Qp$ after possibly twisting by a character. If $s \in \calS^{\st}_*$, then the representation $V(s)$ becomes semistable (non crystalline)  over an abelian extension of $\Qp$ after possibly twisting by a character. If $s \in \calS^{\ngeo}_*$, then $V(s)$ is not a twist of a de Rham representation. In the cases where $V(s)$ is crystalline or semistable, Colmez has explicitly determined in \S 4.5 and \S 4.6 of \cite{PC08} the filtered $\phi$- and $(\phi,N)$-modules $\dcris(V(s))$ and $\dst(V(s))$ in terms of $s$. 

Let us give as an example the description of the parameter $s$ corresponding to the representation $V_p f$ arising from a modular eigenform of level $N$ prime to $p$, weight $k \geq 2$ and character $\eps$. Let $a_p \in \OO_E$ be the eigenvalue of the Hecke operator $T_p$. We assume that $a_p \in \MM_E$ so that $V_p f$ restricted to $D_p$ is irreducible. If $y \in E^\times$, let $\mu_y : \Qp^\times  \to E^\times$ be the character defined by $\mu_y(z)=y^{\vp(z)}$. Let $x_0 : \Qp^\times  \to E^\times$ be the character defined by $x_0(z) = z|z|_p$, so that $x_0(p)=1$ and $x_0(z)=z$ if $z \in \Zp^\times$. The result below then follows from the computations of \S 4.5 of \cite{PC08}.

\begin{theo}\label{tricry}
We have $(V_p f)^* = V(\mu_y,\mu_{\eps(p)/y} x_0^{1-k},\infty)$ where $y \in \MM_E$ is such that $a_p = y + \eps(p) p^{k-1}/y$.
\end{theo}

The equation for $y$ has (in general) two solutions, giving two different parameters $s$ and $s'$ for the same representation. This corresponds to the phenomenon described at the end of theorem \ref{dsirr}.

The constructions of this paragraph have been generalized to $2$-dimensional trianguline representations of $\galk$ by Nakamura in \cite{KN09}, for $K$ a finite extension of $\Qp$.

\section{Arithmetic applications}
\label{aach}

In this chapter, we explain the role that trianguline representations play in the $p$-adic local Langlands correspondence and then in the theory of overconvergent modular forms.

\subsection{The $p$-adic local Langlands correspondence}\label{llsec}

We only give a cursory description of the $p$-adic local Langlands correspondence, and refer to the Bourbaki seminar \cite{LBNB} for a detailed survey and adequate references.

The $p$-adic local Langlands correspondence for $\GL_2(\Qp)$ is a bijection, between certain $2$-dimensional $p$-adic representations of $\galp$, and certain representations of $\GL_2(\Qp)$. The first examples of this correspondence were constructed by Breuil, for semistable and crystalline representations of $\galp$. These examples inspired Colmez to use $(\phi,\Gamma)$-modules in order to give a ``functorial'' construction of these examples, and he realized that the natural condition to impose on the $p$-adic representations which he was considering was that the attached $(\phi,\Gamma)$-module be an extension of two $(\phi,\Gamma)$-modules of rank $1$. This is what led him to define trianguline representations. In the notations of \S\ref{exsec}, if $s \in  \calS_{\irr}$, then the representation of $\GL_2(\Qp)$ corresponding to $V(s)$ by the $p$-adic local Langlands correspondence is a $p$-adic unitary Banach space representation $\Pi(s)$ of $\GL_2(\Qp)$ constructed as follows.

Let $\log_{\Linv}$ be the logarithm normalized by $\log_{\Linv}(p)=\Linv$ (if $\Linv=\infty$, we set $\log_\infty=\vp$) and if $s \in \calS$, let $\delta_s$ be the character $(x|x|_p)^{-1} \delta_1 \delta_2^{-1}$. Note that if $s \in \calS_{\irr}$ then we can have $\Linv \neq \infty$ only if $\delta_s$ is of the form $x^i$ with $i \geq 0$. We can define the notion of a class $\calC^u$ function for $u \in \RR_{\geq 0}$ generalizing the usual case $u \in \ZZ_{\geq 0}$. We denote by $\B(s)$ the space of functions $f : \Qp \to E$ which are of class $\calC^{u(s)}$ and such that $x \mapsto \delta_s(x) f(1/x)$ extends at $0$ to a function of class $\calC^{u(s)}$. The space $\B(s)$ is then endowed with an action of $\GL_2(\Qp)$ given by the formula:
\[ \left[ \pmat{a & b \\ c & d} \cdot f \right] (y) = (x|x|_p\delta_1^{-1})(ad-bc) \cdot \delta_s(cy+d) \cdot f\left( \frac{ay+b}{cy+d}\right). \]

The space $M(s)$ is defined by
\begin{enumerate}
\item if $\delta_s$ is not of the form $x^i$ with $i \geq 0$, then $M(s)$ is the space generated by $1$ and by the functions $y \mapsto \delta_s(y-a)$ with $a \in \Qp$;
\item if $\delta_s$ is of the form $x^i$ with $i \geq 0$, then $M(s)$ is the intersection of $\B(s)$ with the space generated by the functions $y \mapsto \delta_s(y-a)$ and $y \mapsto \delta_s(y-a)\log_{\Linv}(y-a)$ with $a \in \Qp$.
\end{enumerate}
We finally set $\Pi(s)=\B(s)/\widehat{M}(s)$ where $\widehat{M}(s)$ is the closure of $M(s)$ inside $\B(s)$. 

\begin{theo}\label{coltri}
The unitary Banach space representation $\Pi(s)$ of $\GL_2(\Qp)$ is nonzero, topologically irreducible and admissible in the sense of Schneider-Teitelbaum.
\end{theo}

These representations $\Pi(s)$ are called the ``unitary principal series'' and the above theorem is theorem 0.4 of \cite{PC10}. Colmez then proceeds in \cite{PCGL} to attach to any $2$-dimensional $p$-adic representation of $\galp$ a representation of $\GL_2(\Qp)$, and he proves that they have the required properties by using the fact that this is true for trianguline representations, that his construction is suitably continuous, and that trianguline representations are Zariski dense in the deformation space of all $2$-dimensional $p$-adic representations. 

\subsection{Families of Galois representations}\label{locsec}

In this paragraph, we recall the existence of certain rigid analytic spaces which parameterize some families of $p$-adic Galois representations. Recall that the rigid analytic space attached to $\Qp \otimes_{\Zp} \Zp\dcroc{X}$ is the $p$-adic open unit disk and that more generally, the rigid analytic space attached to $E \otimes_{\OO_E} \OO_E \dcroc{X_1,\hdots,X_n}$ is the $n$-dimensional ball over $E$.

We start with a simple example; the group $1+p\Zp$ is topologically generated by the element $1+p$ (unless $p=2$; the constructions of this paragraph can easily be adapted to work when $p=2$). This implies that a character on $1+p \Zp$ is determined by its value at $1+p$. Consider the ring $R=\Zp\dcroc{X}$ and the character $\eta_R : 1+p \Zp \to R^\times$ given by $\eta_R(1+p) = 1+X$. Any character $\eta = 1+p \Zp \to 1+\MM_E$ is obtained from $\eta_R$ by the formula $\eta(g) = f \circ \eta_R(g)$, where $f : R \to \MM_E$ is given by $f(X) = \eta(1+p)-1$. There is therefore a bijection beween the $E$-valued points of the space attached to $\Qp \otimes_{\Zp} \Zp\dcroc{X}$ and the set of characters $\eta : 1+p \Zp \to 1+\MM_E$. The ring $R$ is an example of a \emph{universal deformation ring}, and the rigid analytic space attached to $R[1/p]$ (the $p$-adic open unit disk) parameterizes the family of all $\Qpbar$-valued characters of $1+p\Zp$.

Suppose now that $\overline{\eta} : \Zp^\times \to k_E^\times$ is a character, and let $E_0$ be the smallest extension of $\Qp$ whose residue field is $k_E$ (that is, $E_0 = E \cap \Qpnr$). The natural parameter space for characters $\eta : \Zp^\times \to \Zpbar^\times$ whose reduction modulo $\MM_{\Zpbar}$ is $\overline{\eta}$ is, as above, the rigid analytic space attached to $E_0 \otimes_{\OO_{E_0}} \OO_{E_0} \dcroc{X}$. We denote this space by $\scrX_{\overline{\eta}}$. 

There is likewise a parameter space $\scrX^u_{\overline{\delta}}$ for characters $\delta : \Qp^\times  \to \Qpbar^\times$ which have a fixed slope $u$ and such that $\overline{\delta(p)/p^u} \in k_E^\times$ and $\overline{\delta} \mid_{\Zp^\times} \in k_E^\times$ are fixed, and this parameter space is the rigid analytic space attached to $E_0 \otimes_{\OO_{E_0}} \OO_{E_0} \dcroc{X_1,X_2}$. Denote by $\delta(x)$ the character corresponding to a point $x \in \scrX^u_{\overline{\delta}}$. Colmez proves in \S 5.1 of \cite{PC08} that the representations $V(s)$ live in analytic families of trianguline representations, and his construction has been completed by Chenevier (see \S 3 of \cite{GC10}). 

\begin{theo}\label{trifam}
If $(\delta_1,\delta_2,\infty) \in \calS_{\irr}$, then there exists a neighborhood $\scrU$ of $(\delta_1,\delta_2) \in \scrX^{u_1}_{\overline{\delta}_1} \times \scrX^{u_2}_{\overline{\delta}_2}$ and a free $\OO_{\scrU}$-module $V$ of rank $2$ with an action of $\galp$ such that $V(u)=V(\delta_1(u),\delta_2(u),\infty)$ if $u \in \scrU$.
\end{theo}

Recall that Mazur generalized the construction of $\scrX_{\overline{\eta}}$ in \cite{BM89}, where he proved that for certain groups $G$ and representations $\rhobar : G \to \GL_d(\Fpbar)$, there exists a parameter space $\scrX_{\rhobar}$ for the set of all isomorphism classes of representations $\rho : G \to \GL_d(\Zpbar)$ having reduction modulo $\MM_{\Zpbar}$ isomorphic to $\rhobar$. This applies for example if $\End(\rhobar)=\Fpbar$ and if either $G=\Gal(\QQ_S/\QQ)$ is the Galois group of the maximal extension $\QQ_S$ of $\QQ$ which is unramified outside of a finite set of places $S$, or if $G=\galp$. 

If $G=\galp$ and $d=2$, then for most representations $\rhobar : G \to \GL_d(k_E)$, the rigid analytic space $\scrX_{\rhobar}$ is the one attached to $E_0 \otimes_{\OO_{E_0}} \OO_{E_0} \dcroc{X_1,X_2,X_3,X_4,X_5}$. Theorem \ref{trifam} then implies that inside the $5$-dimensional space $\scrX_{\rhobar}$, there is a countable number of $4$-dimensional subspaces corresponding to trianguline representations. The ``trianguline locus'' of $\scrX_{\rhobar}$ is Zariski dense (it is however a ``thin subset'' of $\scrX_{\rhobar}$ in the terminology of \S 4 of \cite{BC10}). This can be compared with the following result (theorems B and C of \cite{BC08}).

\begin{theo}\label{clsht}
If $b \geq a$, then the locus of $\scrX_{\rhobar}$ corresponding to crystalline (or semistable or de Rham or Hodge-Tate) representations, with Hodge-Tate weights in the range $[a,b]$, is a closed analytic subspace of $\scrX_{\rhobar}$.
\end{theo}

\subsection{Overconvergent modular forms}\label{eksec}

Overconvergent modular forms are objects defined by Coleman in \cite{RC96}, which are $p$-adic generalizations of classical modular forms (for a survey about overconvergent modular forms, see the Bourbaki seminar \cite{ME09}). An ``overconvergent modular form of finite slope'' has a $q$-expansion, which is a $p$-adic limit of $q$-expansions of classical modular eigenforms. One can attach Galois representations to these objects, by taking the limit of the Galois representations attached to the eigenforms in the converging sequence. In this paragraph, we directly define some $p$-adic representations of $\Gal(\Qbar/\QQ)$ by a $p$-adic interpolation process, and merely recall that these representations are the ones which are attached to ``overconvergent modular eigenforms of finite slope''.

Let $N \geq 1$ be an integer prime to $p$ and let $S$ be the set of primes dividing $pN$ and $\infty$. Fix some $2$-dimensional $\Fpbar$-representation $\rhobar$ of $\Gal(\QQ_S/\QQ)$. Let $\scrX^S_{\rhobar}$ be the rigid analytic space attached to the universal deformation space for $\rhobar$, so that every $x \in \scrX^S_{\rhobar}(E)$ corresponds to an $E$-linear representation $V_x$ of $\Gal(\Qbar/\QQ)$, which is unramified outside of $S$, and whose reduction modulo $\MM_E$ is isomorphic to $\rhobar$. For most representations $\rhobar$, $\scrX^S_{\rhobar}$ is a $3$-dimensional rigid analytic ball by results of Weston (see theorem 1 of \cite{TW04}).

Let $\calC_{\mathrm{cl}}$ be the set of points $(x,\lambda) \in \scrX^S_{\rhobar} \times \Gm$ such that $V_x$ is the representation attached to a modular eigenform $f$ on $\Gamma_1(N)$ and $\lambda$ is a root of $X^2-a_p X+\eps(p)p^{k-1}$ (where $k$ is the weight of $f$ and $T_p(f)=a_p f$). Let $\calC$ be the Zariski closure of $\calC_{\mathrm{cl}}$ inside $\scrX^S_{\rhobar} \times \Gm$. By \S 1.5 of \cite{CM98}, we have the following result.

\begin{theo}\label{cmec}
The variety $\calC$ is a rigid analytic curve.
\end{theo}

Coleman and Mazur then show in \cite{CM98} that the Galois representations $V_x$ corresponding to points $(x,\lambda) \in \calC(E)$ are the ones which are attached to the ``level $N$ overconvergent modular eigenforms of finite slope'' defined by Coleman. The curve $\calC$ is called the \emph{eigencurve} (note that the construction of the eigencurve is given in \cite{CM98} assuming that $N=1$ and that $p>2$. The general case is treated in \cite{KB07}). 

The projection of $\calC$ on $\scrX^S_{\rhobar}$ is then a complicated space (for instance, it has infinitely many double points) which is the ``infinite fern'' of \cite{BM97} and \cite{GM98}, see \S 2.5 of \cite{ME09}. The following result (a consequence of theorem 6.3 of \cite{MK03} combined with theorem \ref{dcnz}) describes the restriction to $\galp$ of the representations of $\Gal(\Qbar/\QQ)$ which are constructed in this way.

\begin{theo}\label{kmtr}
If $f$ is an overconvergent modular eigenform of finite slope of level $N$ (i.e.\ if $(V_p f,\lambda) \in \calC(E)$ by the above remark), then $V_p f$ is a trianguline representation.
\end{theo}

The idea is that this theorem is true if $f$ is a classical modular eigenform, by using Saito's theorem, and Kisin deduces that $V_p f$ satisfies the hypothesis of theorem \ref{dcnz} from the classical case by a $p$-adic interpolation argument. 

We then have the following result, Emerton's generalization of the Fontaine-Mazur conjecture for modular forms.

\begin{theo}\label{fmtri}
If $V$ is an irreducible $2$-dimensional $p$-adic representation of $\Gal(\Qbar/\QQ)$, such that
\begin{enumerate}
\item the restriction of $V$ to $I_\ell$ is trivial for all $\ell$ except a finite number,
\item the restriction of $V$ to  $D_p$ is trianguline,
\item $\overline{V}$ satisfies hypotheses (1), (2) and (3') of \S \ref{ptsec},
\end{enumerate}
then $V$ is a twist of the Galois representation attached to an overconvergent cuspidal eigenform of finite slope.
\end{theo}

We now describe the parameter $s \in \calS$ such that $(V_p f)^* = V(s)$, just as we did for classical modular forms at the end of \S \ref{exsec}. Let $f$ be a finite slope overconvergent modular eigenform of level $N$ and character $\eps$. Let $k=w(\det(V_p f))+1 \in E$ (so that if $f$ is classical, then $k$ is the weight of $f$), let $\lambda \in E$ be such that $(V_pf,\lambda) \in \calC$, and let $\mu_\lambda : \Qp^\times \to E^\times$ be the character $z \mapsto \lambda^{\vp(z)}$. The following result is proposition 5.2 of \cite{GC08} (it is a direct consequence of theorem 6.3 of \cite{MK03} and theorem 0.8 of \cite{PC08}). Note that if $k \in \ZZ_{\geq 1}$ and either $\vp(\lambda) =0$ or $\vp(\lambda) =k-1$, then $V_p f$ is reducible in an obvious way.

\begin{theo}\label{gcoe}
We have $(V_p f)^* = V(\delta_1,\det(V_p f)^{-1} \cdot \delta_1^{-1},\Linv_f)$ where
\begin{enumerate}
\item if $k \in \ZZ_{\geq 1}$ and $0 < \vp(\lambda)<k-1$, then $\delta_1=\mu_\lambda$ and if $\Linv_f \neq \infty$, then $V_pf$ is semistable and $\Linv_f$ is the $\Linv$-invariant of $f$;
\item if $k \in \ZZ_{\geq 1}$ and $\vp(\lambda)>k-1$, then $\delta_1=x^{1-k}\mu_\lambda$  and $\Linv_f = \infty$;
\item if $k \notin \ZZ_{\geq 1}$, then $\delta_1=\mu_\lambda$ and $\Linv_f = \infty$.
\end{enumerate}
\end{theo}

Note that case (1) corresponds to $\calS_*^{\cris} \sqcup \calS_*^{\st}$ while cases (2) and (3) correspond to $\calS_*^{\ngeo}$. Coleman's ``small slope criterion'' for the classicality of overconvergent modular eigenforms (\S 6 of \cite{RC96}) can then be interpreted as follows in terms of Galois representations:  if $k \geq 1$ and $0 < \vp(\lambda)<k-1$, then the representation $V_p f$ is potentially semistable at $p$, and therefore the overconvergent modular form $f$ is classical, as predicted by the Fontaine-Mazur conjecture (theorem \ref{emkis}).

We finish this paragraph by discussing the weight-characters of overconvergent cuspidal eigenforms of finite slope. Let $\scrW$ be the \emph{weight space}, that is the parameter space for characters of $\Zp^\times$. The space $\scrW$ is the union of the $p-1$ balls $\scrX_{\overline{\chi}_{\cycl}^i}$ where $0 \leq i \leq p-2$ (unless $p=2$; $\scrW$ is then the union of two balls). If $V$ is a $p$-adic representation of $\galp$, then by class field theory $x_0^{-1} \cdot \det(V)$ gives a character of $\Qp^\times$ whose restriction to $\Zp^\times$ is the \emph{weight-character} $\kappa_V$ of $V$. This gives rise to a map $\kappa : \scrX^S_{\rhobar} \to \scrW$ and by composition to a map $\calC \to \scrW$ which satisfies the following property by \S 1.5 of \cite{CM98}.

\begin{theo}\label{egwt}
The map $\calC \to \scrW$ is, locally in the domain, finite and flat.
\end{theo}

We now explain that if $N=1$ and  $p \in \{2,3,5,7\}$, then one can give a ``local'' realization of the eigencurve. A point $(\kappa,\lambda) \in \scrW \times \Gm$ is said to be \emph{special} if $\kappa = x_0^k$ for some $k \geq 2$ and $\lambda^2=p^{k-2}$. Let $\scrW \widetilde{\times} \Gm$ be the blow-up of $\scrW \times \Gm$ at the special points (so that $\scrW \widetilde{\times} \Gm$ can be seen as a subspace of the space $\calS$ of \S\ref{exsec}). 

Consider the map $\calC \to \scrW \widetilde{\times} \Gm$ given by $(V_x,\lambda) \mapsto (\kappa_x,\lambda,\Linv_x)$ at the special points and by $(V_x,\lambda) \mapsto (\kappa_x,\lambda)$ elsewhere. The following theorem is the main result of \cite{GC08}.

\begin{theo}\label{egcol}
The map $\calC \to \scrW \widetilde{\times} \Gm$ is a rigid analytic map, and if $N=1$ and $p \in \{2,3,5,7\}$, then it is a closed immersion.
\end{theo}

Some important ideas underlying the proof of this theorem are Colmez' theorem 0.5 of \cite{PCLI} expressing the $\Linv$-invariant as the derivative of the $U_p$-eigenvalue, the fact that if $p \in \{2,3,5,7\}$ and $S=\{p,\infty\}$, then an odd $2$-dimensional $p$-adic representation of $\Gal(\QQ_S/\QQ)$ is determined by its restriction to $\galp$ (proposition 1.8 of \cite{GC08}), and a local study of families of trianguline representations.

\subsection{Trianguline representations and Selmer groups}\label{bcsec}

Since the $(\phi,\Gamma)$-module attached to a trianguline representation $V$ has a particularly simple structure, one can use this structure to study the cohomology groups attached to $V$, in particular the Selmer group and its variants. Some of the techniques which are available in the ordinary case for that study (such as \cite{RG89}) can be extended to the case of trianguline representations. 

For example, it is possible to give a generalized definition of the usual $\Linv$-invariant (see Benois' \cite{DB09}, where an $\Linv$-invariant is constructed for representations which are not necessarily $2$-dimensional), and to study the Selmer groups corresponding to families of trianguline representations, such as those carried by the eigencurve or more general eigenvarieties (as in the book \cite{BC09} by Bella\"{\i}che and Chenevier, and in Pottharst's \cite{JP08} and \cite{JP10}). In this way, it is possible to prove some new cases of the Bloch-Kato conjectures, by establishing some ``lower semicontinuity'' results about the rank of the Selmer groups (see \cite{BC09} and Bella\"{\i}che's \cite{JB10}). 

The systematic study of families of trianguline representations, in connection with the theory of families of automorphic forms, is an increasingly important topic which we do not say anything more about, because it is rapidly progressing and deserves a survey of its own.

\noindent \textbf{Acknowledgements}: I had several very helpful conversations with Ga\"etan Chenevier. The two referees' detailed comments allowed me to correct many (hopefully, all) inaccuracies.

\providecommand{\bysame}{\leavevmode ---\ }
\providecommand{\og}{``}
\providecommand{\fg}{''}
\providecommand{\smfandname}{and}
\providecommand{\smfedsname}{eds.}
\providecommand{\smfedname}{ed.}
\providecommand{\smfphdthesisname}{PhD}

\end{document}